\documentclass[a4paper,12pt]{amsart}
\usepackage{graphics}
\usepackage{amssymb}
\usepackage{amsmath}
\usepackage{latexsym}
\usepackage{amscd}
\usepackage{mathrsfs}
\newtheorem{Th}{Theorem}
\newtheorem{Prop}[Th]{Proposition}

\newtheorem{Rem}[Th]{\sc Remark}
\newtheorem{Qu}[Th]{\sc Questions}

\newtheorem{Cor}[Th]{Corollary}
\newtheorem{Lemma}[Th]{Lemma}
\newcommand{\be}{\begin{eqnarray*}}
\newcommand{\ee}{\end{eqnarray*}}

\newcommand{\LccEF}{{\mathcal L}_{\mathop{\rm cc}\nolimits}(E,F)}

\newcommand{\Lcc}{{\mathcal L}_{\mathop{\rm cc}\nolimits}(}

\newcommand{\Nuc}{{\mathcal N}(}
\newcommand{\I}{{\mathcal I}(}

\newcommand{\Proof}{\noindent {\sc Proof. }}
\newcommand{\fin}{\hspace*{\fill} $\Box$}
\newcommand{\finesp}{\hspace*{\fill} $\Box$\vspace{.5\baselineskip}}

\newcommand{\ra}{\rightarrow}

\newcommand{\Ra}{\Rightarrow}       
\newcommand{\ptp}{\widehat{\otimes}_\pi}
\newcommand{\pti}{\widehat{\otimes}_\epsilon}
\newcommand{\espv}{\vspace{.5\baselineskip}}

\begin{document}
\title{The Dunford-Pettis property on tensor products}

\author[M. Gonz\'alez]{Manuel Gonz\'alez}
\address{Departamento de Matem\'aticas \\
      Facultad de Ciencias\\
      Universidad de Cantabria \\ 39071 Santander (Spain)}
\email{gonzalem@ccaix3.unican.es}
\thanks{The first named author was supported
in part by DGICYT Grant PB 97--0349 (Spain)}

\author[J. M. Guti\'errez]{Joaqu\'\i n M. Guti\'errez}
\address{Departamento de Matem\'atica Aplicada\\
      ETS de Ingenieros Industriales \\
      Universidad Polit\'ecnica de Madrid\\
      C. Jos\'e Guti\'errez Abascal 2 \\
      28006 Madrid (Spain)}
\email{jgutierrez@etsii.upm.es}
\thanks{The second named author was supported
in part by DGICYT Grant PB 96--0607 (Spain)}
\thanks{\hspace*{\fill}\scriptsize file dpptp.tex}




\begin{abstract}
We show that, in some cases, the projective and the injective
tensor products of two Banach spaces do not have the
Dunford-Pettis property (DPP). As a consequence, we obtain  that
$(c_0\ptp c_0)^{**}$ fails the DPP. Since $(c_0\ptp c_0)^{*}$ does
enjoy it, this provides a new space with the DPP whose dual fails
to have it. We also prove that, if $E$ and $F$ are ${\mathscr
L}_1$-spaces, then $E\pti F$ has the DPP if and only if both $E$
and $F$ have the Schur property. Other results and examples are
given.
\end{abstract}

\maketitle

A Banach space $E$ has the {\it Dunford-Pettis property\/} ({\it
DPP,} for short) if every weakly compact operator on $E$ is {\it
completely continuous,} i.e., takes weak Cauchy sequences into
norm Cauchy sequences \cite{Gr}. Equivalently, $E$ has the DPP if,
for all weakly null sequences $(x_n)\subset E$ and
$(\phi_n)\subset E^*$ (the dual of $E$), we have
$\lim\phi_n(x_n)=0$. If $E^*$ has the DPP, then so does $E$, but
the converse is not true \cite{S}. A Banach space with the Schur
property has the DPP. The DPP is inherited by complemented
subspaces. For more on the DPP, the reader is referred to
\cite{DiDP,CG}.

It was unknown for many years if the Dunford-Pettis property of
$E$ and $F$ implies that of their projective tensor product $E\ptp
F$ and of their injective tensor product $E\pti F$ \cite{DiDP}.
Talagrand \cite{T} found a Banach space $E$ so that $E^*$ has the
Schur property but $C([0,1],E)=C[0,1]\pti E$ and
$L^1([0,1],E^*)=L^1[0,1]\ptp E^*$ fail the DPP.

It is proved in \cite{Ry} that, if $E$ and $F$ have the DPP and
contain no copy of $\ell_1$, then $E\ptp F$ has both properties.
It is shown in \cite{LP} that, if $E$ and $F$ have the Schur
property, then $E\pti F$ has the Schur property.

Answering a question of \cite{CG}, it is proved in \cite{BV} that,
for compact spaces $K_1,\ldots,K_n$, the space
$C(K_1)\ptp\cdots\ptp C(K_n)$ has the DPP if and only if
$K_1,\ldots,K_n$ are scattered.

Taking advantage of an idea of \cite{BV} we prove, among others,
the following results, where ${\mathcal L}(E,F)$ denotes the space
of all (linear bounded) operators from $E$ into $F$ while $\LccEF$
is the space of all completely continuous operators from $E$ into
$F$:

(a) Suppose $E$ is not Schur, $F$ contains a copy of $\ell_1$, and
${\mathcal L}(E,F^*)=\Lcc E,F^*)$. Then $E\ptp F$ does not have
the DPP.

(b) Suppose $E^*$ does not have the Schur property, $F^*$ contains
a copy of $\ell_1$, and ${\mathcal L}(E^*,F^{**})=\Lcc
E^*,F^{**})$. Then $(E\pti F)^*$ does not have the DPP. As a
consequence, $(c_0\ptp c_0)^{**}$ fails the DPP, which answers a
question of
\cite{CG}.

(c) If $E$ and $F^*$ do not have the Schur property, and
${\mathcal L}(F^*,E^{**})=\Lcc F^*,E^{**})$, then $E\pti F$ fails
the DPP. As a consequence, if $E, F$ are ${\mathscr L}_1$-spaces,
then $E\pti F$ has the Dunford-Pettis property if and only if it
has the Schur property, if and only if both $E$ and $F$ have the
Schur property.

(d) Assume $E^*$ is not Schur, $E^{**}$ contains no complemented
copy of $\ell_1$, and $F^*$ contains a complemented copy of
$\ell_1$. Then $(E\ptp F)^*$ does not have the DPP.\espv

Throughout, $E$ and $F$ will denote Banach spaces, $e_n$ denotes,
as usual, the vector $(0,\ldots,0,1,0\ldots)$ with $1$ in the
$n\,$th place. We shall use the following notation for some
subspaces of ${\mathcal L}(E,F)$: ${\mathcal K}(E,F)$ for the
space of compact operators, $\Nuc E,F)$ for the nuclear operators
and $\I E,F)$ for the integral operators. The nuclear norm of an
operator $T$ is denoted by $\| T\|_{\mathop{\rm nuc}\nolimits}$,
while $\| T\|_{\mathop{\rm int}\nolimits}$ stands for the integral
norm. Given an operator $T\in{\mathcal L}(E,F)$, its adjoint is
denoted by $T^*\in{\mathcal L}(F^*,E^*)$. By $E\simeq F$ we mean
that $E$ and $F$ are isomorphic.\espv

The {\it Banach-Mazur distance $d(E,F)$} between two isomorphic
Banach spaces $E$ and $F$ is defined by $\inf
\left(\|T\|\|T^{-1}\|\right)$ where the infimum is taken over all
isomorphisms $T$ from $E$ onto $F$. Recall that a Banach space $E$
is an {\it ${\mathscr L}_\infty$-space\/} (resp.\ {\it ${\mathscr
L}_1$-space\/})
\cite{B} if there is $\lambda\geq 1$ such that every finite
dimensional subspace of $E$ is contained in another subspace $N$
with $d\left( N,\ell_\infty^n\right)\leq\lambda$ (resp.\ $d\left(
N,\ell_1^n\right)\leq\lambda$) for some integer $n$. The
${\mathscr L}_\infty$-spaces, the ${\mathscr L}_1$-spaces and all
their dual spaces have the DPP
\cite[Corollary 1.30]{B}.

Since some of our results and examples are on ${\mathscr
L}_\infty$ and ${\mathscr L}_1$-spaces, we shall first recall a
well known result on tensor products of these spaces. Observe that
it provides examples of tensor products with the DPP.

\begin{Prop}{\rm \cite[Theorem~34.9]{DF}}

{\rm (a)} If $E$ and $F$ are ${\mathscr L}_\infty$-spaces, then
$E\pti F$ is an ${\mathscr L}_\infty$-space.

{\rm (b)} If $E$ and $F$ are ${\mathscr L}_1$-spaces, then $E\ptp
F$ is an ${\mathscr L}_1$-space.
\end{Prop}

It is also known that, if $E$ is an ${\mathscr L}_1$-space and $F$
is an ${\mathscr L}_\infty$-space, then $E\pti F$ and $E\ptp F$
have the DPP (see \cite{Bo2}, \cite{Em1} and \cite{Em2}).

 From the representation $(E\ptp F)^*= {\mathcal L}(E,F^*)$, we
obtain:

\begin{Lemma}
\label{wnull}{\rm \cite{BV}}
Assume ${\mathcal L}(E,F^*)=\Lcc E,F^*)$, and take a weakly null
sequence $(x_n)\subset E$ and a bounded sequence $(y_n)\subset F$.
Then the sequence $(x_n\otimes y_n)$ is weakly null in $E\ptp F$.
\end{Lemma}

The following result is essentially contained in \cite{BV}:

\begin{Th}
\label{proy}
Suppose $E$ does not have the Schur property, $F$ contains a copy
of $\ell_1$, and ${\mathcal L}( E,F^*)=\Lcc E,F^*)$. Then $E\ptp
F$ does not have the DPP.
\end{Th}

\Proof
Since $E$ does not have the Schur property, we can find a
normalized weakly null sequence $(x_n)\subset E$ that may be
supposed to be basic. Let $(\phi_n)\subset E^*$ be a bounded
sequence with $\phi_i(x_j)=\delta_{ij}$.

Since $F$ contains a copy of $\ell_1$, there is a surjective
operator $q:F\ra\ell_2$ \cite[Proposition~3]{OP}. Consider the
operator $T:E\ptp F\ra\ell_2$ given by
$$
T(x\otimes y):=\left( \phi_k(x)\langle
q(y),e_k\rangle\right)_{k=1}^\infty .
$$
Since
$$
\| T(x\otimes y)\|\leq \sup_k\|\phi_k\|\cdot\|
q\|\cdot\|x\|\cdot\|y\|\, ,
$$
$T$ is a well defined weakly compact operator. Choose a bounded
sequence $(y_n)\subset F$ so that $q(y_n)=e_n$. Then
$$
T(x_n\otimes y_n)=\left( \phi_k(x_n)\langle
q(y_n),e_k\rangle\right)_{k=1}^\infty =e_n .
$$
So, by Lemma~\ref{wnull}, $T$ is not completely continuous, and
$E\ptp F$ fails the DPP.\finesp

If we replace the condition ``$F$ contains a copy of $\ell_1$'' by
the more general ``$F$ has a quotient isomorphic to $\ell_p$
$(1<p<\infty)$'', the proof still works. However, the Theorem does
not gain in extension. Indeed, we can assume that $F$ has the DPP.
Then, if $F$ has a quotient isomorphic to $\ell_p$ $(1<p<\infty)$,
it must contain a copy of $\ell_1$. This is also true for
Theorems~\ref{newproy} and \ref{ptidual}.

Theorem~\ref{proy} applies, for instance, if $E$ is an ${\mathscr
L}_\infty$-space without the Schur property and $F$ is an
${\mathscr L}_\infty$-space containing a copy of $\ell_1$: for the
equality ${\mathcal L}( E,F^*)=\Lcc E,F^*)$, see
\cite[Theorem~3.7]{DJT}.

We now give another interesting example. Let ${\mathbb T}$ be the
unit circle with its normalized Lebesgue measure; let $H^1$ denote
the closed subspace of $L^1=L^1({\mathbb T})$ spanned by the
functions $\{ e^{int}:n\geq 0\}$. Although $L^1/H^1$ contains an
isomorphic copy of $L^1$, every operator $L^1/H^1\ra L^1$ is
completely continuous \cite{Bo3}. Let $E=L^1/H^1$ and $F=L^\infty
[0,1]$. We know that $F^*$ may be written as $L^1(\mu)$ for some
measure $\mu$. Given an operator $T:E\ra F^*\equiv L^1(\mu)$,
since $\overline{T(E)}$ is separable, $\overline{T(E)}$ is
isomorphic to a subspace of $L^1[0,1]$
\cite[Theorem~IV.1.7]{GD}. Therefore, $T$ is completely
continuous. Applying Theorem~\ref{proy}, we have that $L^1/H^1\ptp
L^\infty [0,1]$ does not have the DPP, although $L^1[0,1]\ptp
L^\infty [0,1]$ does have it \cite{Bo2}. Note also that
$L^1/H^1\ptp L^1(\mu)$ has the DPP for all $\mu$ \cite{Ra}.

\begin{Rem}
{\rm
Theorem~\ref{proy} does not include Talagrand's example \cite{T}
of a space $E$ with the Schur property such that
$L^1([0,1],E)=L^1[0,1]\ptp E$ fails the DPP.
}
\end{Rem}
Indeed,
Theorem~\ref{proy} cannot be applied to the space $L^1[0,1]$,
since we have the following facts:

(a) If ${\mathcal L}(L^1[0,1],F^*)=\Lcc L^1[0,1],F^*)$, then $F$
contains no copy of $\ell_1$.

(b) If ${\mathcal L}(E,L^1[0,1]^*)=\Lcc E,L^1[0,1]^*)$, then $E$
has the Schur property.

For (a), suppose $F$ has a subspace $M\simeq \ell_1$. Then
$F^*/M^\perp\simeq\ell_\infty$ contains a copy of $L^1[0,1]$. If
$T:L^1[0,1]\to F^*/M^\perp$ is an isomorphism, then we can write
$T=Q\circ S$ \cite[Proposition~1]{GrL1}, where $Q:F^*\to
F^*/M^\perp$ is the quotient map. Therefore, $S:L^1[0,1]\to F^*$
is not completely continuous.

To prove (b), take a normalized weakly null, basic sequence
$(x_n)\subset E$, and a bounded sequence $(\phi_n)\subset E^*$
with $\phi_i(x_j)=\delta_{ij}$. Since $L^1[0,1]^*$ is isomorphic
to $\ell_\infty$, the mapping
$$
x\in E\longmapsto (\phi_n(x))_{n=1}^\infty\in\ell_\infty
$$
provides an operator $E\to L^1[0,1]^*$ which is not completely
continuous.

\begin{Cor}
Assume $E, F$ are infinite dimensional ${\mathscr
L}_\infty$-spaces, and $E\ptp F$ has the DPP. Then, either $E$ and
$F$ have the Schur property or $E^*$ and $F^*$ have the Schur
property.
\end{Cor}

\Proof
Observe first that ${\mathcal L}(E,F^*)=\Lcc E,F^*)$
\cite[Theorem~3.7]{DJT}.
If $E$ does not have the Schur property, then, by
Theorem~\ref{proy}, $F$ contains no copy of $\ell_1$, so $F^*$ has
the Schur property. Hence, $F$ fails the Schur property and, by
the same argument, $E^*$ has the Schur property. On the other
hand, if $E$ has the Schur property, then it contains a copy of
$\ell_1$. Therefore, by Theorem~\ref{proy}, $F$ has the Schur
property.\fin

\begin{Rem}{\rm
If $E^*$ and $F^*$ have the Schur property, then so does $(E\ptp
F)^*$ \cite{Ry}. Hence, $E\ptp F$ has the DPP. Assume now that $E$
and $F$ have the Schur property. We do not know if this implies
that $E\ptp F$ has the Schur property. Consider, for instance, an
${\mathscr L}_\infty$-space $X$ with the Schur property \cite{BD}.
We do not know if $X\ptp X$ has the Schur property. We shall see,
however, that $(X\ptp X)^*$ fails the DPP (see
Questions~\ref{qq}).}
\end{Rem}

\begin{Cor}
\label{BoDe}
Let $E, F$ be infinite dimensional separable ${\mathscr
L}_\infty$-spaces. Then, the space $(E\ptp F)^*$ has the DPP if
and only if $E^*\simeq F^*\simeq\ell_1$.
\end{Cor}

\Proof
If $E^*\simeq F^*\simeq\ell_1$, then $E$ contains no copy of
$\ell_1$, so $(E\ptp F)^*={\mathcal L}(E,F^*)$ has the Schur
property \cite[Theorem~3.3]{Ry}. Conversely, if
$E^*\not\simeq\ell_1$ and $E$ is infinite dimensional and
separable, then $E^*\simeq C[0,1]^*$ \cite[Theorem 3.1]{B}.
Hence,
$$
(E\ptp F)^*\simeq {\mathcal L}(F,C[0,1]^*)=(F\ptp C[0,1])^*.
$$
Now, Theorem~\ref{proy} implies that $F\ptp C[0,1]$ does not have
the DPP. So, neither $(E\ptp F)^*$ does.\finesp

We shall now give a variant of Theorem~\ref{proy} for which we
need a lemma. Recall that $E$ contains no complemented copy of
$\ell_1$ if and only if ${\mathcal L}(E,\ell_1)={\mathcal
K}(E,\ell_1)$ \cite[Exercise~VII.3 and Theorem~V.10]{Di}.

\begin{Lemma}
Assume $E$ contains no complemented copy of $\ell_1$, and $F$
contains a sequence $(y_n)\subset F$ equivalent to the
$c_0$-basis. Then, for every  bounded sequence $(x_n)\subset E$,
the sequence $(x_n\otimes y_n)$ is weakly null in $E\ptp F$.
\end{Lemma}

\Proof
Let $M$ be the closed subspace generated by $\{ y_n\}$. Take $T\in
{\mathcal L}(E, M^*)={\mathcal K}(E, M^*)$. Since $(y_n)$ is
weakly null, we have $\langle Tx_n,y_n\rangle\ra 0$ and so
$(x_n\otimes y_n)$ is weakly null in $E\ptp M$. The linear mapping
$E\ptp M\to E\ptp F$ given by $x\otimes y\mapsto x\otimes y$ is
continuous. Therefore, $(x_n\otimes y_n)$ is weakly null in $E\ptp
F$.\fin

\begin{Th}
\label{newproy}
Suppose $E$ contains a copy of $\ell_1$, but $\ell_1$ is not
complemented in $E$, and $F$ contains a copy of $c_0$. Then $E\ptp
F$ does not have the DPP.
\end{Th}

\Proof
Consider a surjective operator $q:E\ra\ell_2$ and select a bounded
sequence $(x_n)\subset E$ with $q(x_n)=e_n$. Let $(y_n)\subset F$
be a sequence equivalent to the $c_0$-basis, and take a bounded
sequence $(\psi_n)\subset F^*$ with $\psi_i(y_j)=\delta_{ij}$.
Define $T:E\ptp F\ra\ell_2$ by
$$
T(x\otimes y):=\left(\psi_k(y)\langle
q(x),e_k\rangle\right)_{k=1}^\infty .
$$
The proof proceeds as in Theorem~\ref{proy}.\finesp

As a consequence of Theorem~\ref{newproy}, letting
$E=\ell_\infty$, $F=c_0\pti L^1[0,1]$, the space $E\ptp F$ fails
the DPP. Theorem~\ref{proy} is not applicable to this example. On
the other hand, if $E$ is a somewhat reflexive ${\mathscr
L}_\infty$-space
\cite{BD} and $F=\ell_\infty$, then Theorem~\ref{newproy} is not
applicable, but Theorem~\ref{proy} implies that $E\ptp F$ fails
the DPP.

\begin{Th}
\label{ptidual}
Suppose $E^*$ does not have the Schur property, $F^*$ contains a
copy of $\ell_1$, and ${\mathcal L}(E^*,F^{**})=\Lcc E^*,F^{**})$.
Then $(E\pti F)^*$ does not have the DPP.
\end{Th}

\Proof
Since $E^*$ is not Schur, we can find a weakly null normalized,
basic sequence $(\phi_n)\subset E^*$. Let $(z_n)\subset E^{**}$ be
a sequence with $\|z_n\|\leq C$ and $\langle \phi_i,z_j\rangle
=\delta_{ij}$. Since $F^*$ contains a copy of $\ell_1$, there is a
surjective operator $q:F^*\ra\ell_2$.

Recall that $(E\pti F)^*=\I E,F^*)$
\cite[Corollary~VIII.2.12]{DU}. Define $T:\I E,F^*)\ra\ell_2$ by
$$
T(A)=\left(\left\langle
A^*q^*(e_i),z_i\right\rangle\right)_{i=1}^\infty .
$$

Since $A^*$ is integral \cite[Corollary~VIII.2.11]{DU}, and
$\ell_2$ has the Radon-Nikod\'ym property and the approximation
property, we have
$$
A^*q^*\in\I \ell_2,E^*)=\Nuc\ell_2,E^*)=\ell_2\ptp E^*,
$$
where the equalities are isometric isomorphisms (see
\cite[Theorem~VIII.4.6]{DU} and \cite[p.~3]{Pi}).

Assume therefore $A^*q^*=\sum_{n=1}^\infty y_n\otimes \xi_n$, with
$y_n=(y_n(i))_{i=1}^\infty\in\ell_2$, $\xi_n\in E^*$ and
$\sum_{n=1}^\infty \| y_n\|\cdot\|\xi_n\|<+\infty$. Then
\be
\|T(A)\|&=&\left\|\left(\left\langle\left(\sum_{n=1}^\infty
y_n\otimes \xi_n\right)
(e_i),z_i\right\rangle\right)_{i=1}^\infty\right\|_2\\
&=&\left\|\left(\sum_{n=1}^\infty
y_n(i)\langle\xi_n,z_i\rangle\right)_{i=1}^\infty\right\|_2\\
&\leq&\sum_{n=1}^\infty\left\|\left(
y_n(i)\langle\xi_n,z_i\rangle\right)_{i=1}^\infty\right\|_2\\
&\leq&C\sum_{n=1}^\infty\left\|y_n\right\|
\cdot\left\|\xi_n\right\| .
\ee
Taking the infimum over all representations of $A^*q^*$ of the
above form, we get
$$
\|T(A)\|\leq C\|A^*q^*\|_{\mathop{\rm nuc}\nolimits}
=C\|A^*q^*\|_{\mathop{\rm int}\nolimits}
\leq C\|A\|_{\mathop{\rm int}\nolimits}\cdot\|q\|
$$
(see \cite[Theorem~VIII.2.7]{DU}), so $T$ is a well defined weakly
compact operator.

Now, take a bounded sequence $(\psi_n)\subset F^*$, with
$q(\psi_n)=e_n$. By Lemma~\ref{wnull}, the sequence
$(\phi_n\otimes\psi_n)$ is weakly null in $E^*\ptp F^*$ and so in
$\I E,F^*)$. However,
\be
T(\phi_n\otimes\psi_n)&=&\left(\left\langle
(\phi_n\otimes\psi_n)^*q^*(e_i),
z_i\right\rangle\right)_{i=1}^\infty\\
&=&\left(\left\langle
q^*(e_i),(\phi_n\otimes\psi_n)^{**}(z_i)
\right\rangle\right)_{i=1}^\infty\\
&=&\left(\left\langle
q^*(e_i),\delta_{in}\psi_n
\right\rangle\right)_{i=1}^\infty\\
&=&\left(\left\langle
e_i,\delta_{in}q(\psi_n)
\right\rangle\right)_{i=1}^\infty\\
&=&\left(\left\langle
e_i,\delta_{in}e_n
\right\rangle\right)_{i=1}^\infty\\
&=&e_n.
\ee
Therefore, $T$ is not completely continuous.\fin

\begin{Cor}
\label{CK}
The space $(C(K_1)\ptp C(K_2))^{**}$ does not have the DPP, for
all infinite compact Hausdorff spaces $K_1$ and $K_2$.
\end{Cor}

\Proof
If one of the compact spaces is not scattered, then $C(K_1)\ptp
C(K_2)$ fails the DPP \cite{BV}. A fortiori, its bidual fails the
DPP.

Assume both compact spaces are scattered. Since $C(K_1)^*$ has the
approximation property \cite[Example~VIII.3.11]{DU}, we have
\cite[Corollary~5.3]{DF}:
$$
(C(K_1)\ptp C(K_2))^*={\mathcal L}(C(K_1),C(K_2)^*)={\mathcal
K}(C(K_1), C(K_2)^*)= C(K_1)^*\pti C(K_2)^*.
$$
Therefore, by Theorem~\ref{ptidual}, the space $(C(K_1)\ptp
C(K_2))^{**}=(C(K_1)^*\pti C(K_2)^*)^*$ does not have the
DPP.\finesp

Observe that, if $K_1$ and $K_2$ are scattered, $(C(K_1)\ptp
C(K_2))^{*}$ has the Schur property \cite{Ry} and, hence, the DPP.
This gives a new family of spaces with the DPP whose duals fail
the DPP. Up to now, there was essentially one example of such
spaces, namely $\ell_1(\ell_2^n)$, due to Stegall \cite{S}.

Theorem~\ref{ptidual} implies that, if $E, F$ are infinite
dimensional ${\mathscr L}_1$-spaces, then $(E\pti F)^*$ fails the
DPP.
 From our next result we shall obtain that even $E\pti F$
fails the DPP unless $E$ and $F$ have the Schur property.

\begin{Th}
\label{inypr}
Suppose $E$ and $F^*$ do not have the Schur property, and
${\mathcal L}(F^*,E^{**})=\Lcc F^*,E^{**})$. Then $E\pti F$ does
not have the DPP.
\end{Th}

\Proof
Choose a weakly null normalized, basic sequence $(x_n)\subset E$,
and $(\phi_n)\subset E^*$ bounded so that
$\phi_i(x_j)=\delta_{ij}$.

We can find a weakly null normalized, basic sequence
$(\psi_n)\subset F^*$. By \cite[Remark~III.1]{JR}, after passing
to a subsequence, there is a bounded sequence $(f_n)\subset
[\psi_n]^*$ with $f_i(\psi_j)=\delta_{ij}$, where $[\psi_n]$ is
the closed linear span of the sequence $(\psi_n)$ such that, when
defining $T:F\ra [\psi_n]^*$ by $(Ty)(x)=x(y)$ for all $y\in F$
and $x\in [\psi_n]$, we have $T(F)\supseteq [f_n]$. The operator
$\gamma: T^{-1}([f_n])/\ker T\ra [f_n]$ given by $\gamma(y+\ker
T)=T(y)$, for all $y\in T^{-1}([f_n])$, is an isomorphism.
Therefore, there is a bounded sequence $(y_n)\subset F$ such that
$T(y_n)=f_n$. Moreover,
$$
\delta_{ij}=f_i(\psi_j)=(Ty_i)(\psi_j)=\psi_j(y_i).
$$

Since every operator in ${\mathcal I}(E,F^*)=(E\pti F)^*$ is
completely continuous \cite[Theorem~VIII.2.9]{DU}, the sequence
$(x_n\otimes y_n)$ is weakly null in $E\pti F$. By
Lemma~\ref{wnull}, the sequence $(\phi_n\otimes\psi_n)$ is weakly
null in $E^*\ptp F^*$. Let $J:E^*\ptp F^*\to (E\pti F)^*$ be the
natural mapping. Then the sequence $\left(
J(\phi_n\otimes\psi_n)\right)$ is weakly null in $(E\pti F)^*$.
However, $\langle J(\phi_n\otimes\psi_n), x_n\otimes y_n\rangle
=1$, so $E\pti F$ fails the DPP.\fin

\begin{Cor}
\label{iny}
Suppose $E$ does not have the Schur property, $E^*$ contains no
complemented copy of $\ell_1$, and $F$ contains a complemented
copy of $\ell_1$. Then $E\pti F$ does not have the DPP.
\end{Cor}

\Proof
Assume first $F=\ell_1$. Since $E^*$ contains no complemented copy
of $\ell_1$, $E^{**}$ contains no copy of $\ell_\infty$.
Therefore, ${\mathcal L}(\ell_\infty,E^{**})={\mathcal
L}_{\mathop{\rm cc}\nolimits} (\ell_\infty,E^{**})$. By
Theorem~\ref{inypr}, $E\pti\ell_1$ does not have the DPP. Since
$E\pti F$ contains a complemented copy of $E\pti\ell_1$, the
result is proved.\fin

\begin{Cor}
\label{Schur}
Let $E$ and $F$ be infinite dimensional ${\mathscr L}_1$-spaces.
The following assertions are equivalent:

{\rm (a)} $E\pti F$ has the DPP;

{\rm (b)} $E\pti F$ has the Schur property;

{\rm (c)} both $E$ and $F$ have the Schur property.
\end{Cor}

\Proof
(a) $\Ra$ (b). Suppose $E\pti F$ does not have the Schur property.
Then  at least one of the spaces, say $E$, fails the Schur
property \cite{LP}. Moreover, $E^*$ contains no complemented copy
of $\ell_1$ \cite[Theorem~III]{LR} and $F$ contains a complemented
copy of $\ell_1$ \cite[Theorem~I]{LR}. By Corollary~\ref{iny}, the
space $E\pti F$ fails the DPP.

(b) $\Ra$ (c) is clear.

(c) $\Ra$ (a). The injective tensor product of two spaces with the
Schur property has the Schur property \cite{LP} and therefore the
DPP.\finesp

Finally, we give a result on the dual of the projective tensor
product.

\begin{Th}
\label{ptpdual}
Assume $E^*$ does not have the Schur property, $E^{**}$ contains
no complemented copy of $\ell_1$, and $F^*$ contains a
complemented copy of $\ell_1$. Then $(E\ptp F)^*$ does not have
the DPP.
\end{Th}

\Proof
We assume first that $F^*\simeq \ell_1$. Take a weakly null
normalized, basic sequence $(\phi_n)\subset E^*$, and
$(z_n)\subset E^{**}$ bounded with $\langle z_i,\phi_j\rangle
=\delta_{ij}$. We have $(E\ptp F)^*\simeq {\mathcal L}(E,\ell_1)$.

By Kalton's test \cite{K}, the sequence $(T_n)=(\phi_n\otimes
e_n)$ is weakly null in ${\mathcal K}(E,\ell_1)$, and so in
${\mathcal L}(E,\ell_1)$.

Since $E^{**}$ contains no complemented copy of $\ell_1$,
$E^{***}$ contains no copy of $\ell_\infty$. Hence, ${\mathcal
L}(\ell_\infty,E^{***})=\Lcc \ell_\infty,E^{***})$. By
Lemma~\ref{wnull}, the sequence $(z_n\otimes e_n)$ is weakly null
in $E^{**}\ptp\ell_\infty$. Consider the operator
$\gamma:E^{**}\ptp\ell_\infty\to {\mathcal L}(E,\ell_1)^*$ given
by
$$
\langle\gamma(z\otimes\xi),T\rangle =\langle T^{**}(z),\xi\rangle
\quad\mbox{for }T\in {\mathcal L}(E,\ell_1), z\in E^{**},
\xi\in\ell_\infty .
$$
Then,
$$
\langle\gamma (z_n\otimes e_n),\phi_n\otimes e_n\rangle =\langle
z_n,\phi_n\rangle\cdot\langle e_n,e_n\rangle =1.
$$
Therefore, ${\mathcal L}(E,\ell_1)$ does not have the DPP.

In general, if $P:F^*\to F^*$ is a projection with
$P(F^*)\simeq\ell_1$, then $Q(T):=P\circ T$ defines a projection
on ${\mathcal L}(E,F^*)$ whose range is ${\mathcal L}(E,P(F^*))$.
Therefore, ${\mathcal L}(E,F^*)$ does not have the DPP.\finesp

\begin{Cor}
\label{Linf}
Let $E$ and $F$ be infinite dimensional ${\mathscr
L}_\infty$-spaces, at least one of which contains a copy of
$\ell_1$. Then $(E\ptp F)^*$ does not have the DPP.
\end{Cor}

\begin{Qu}
\label{qq}
{\rm Consider an infinite dimensional separable ${\mathscr
L}_\infty$-space $X$ with the Schur property such that $X^*$ is
isomorphic to $C[0,1]^*$ \cite{BD}.

By Theorem~\ref{iny}, $X^*\pti X^*$ fails the DPP. By
Corollary~\ref{BoDe}, $(X\ptp X)^*$ fails the DPP. We raise the
following questions: does $X\ptp X$ have the Schur property?, does
$X\ptp X$ have the DPP?}
\end{Qu}


\begin{thebibliography}{99}

\bibitem{BV} F. Bombal and I. Villanueva, On the Dunford-Pettis
property of the tensor product of $C(K)$ spaces, to appear in {\it
Proc.\ Amer.\ Math.\ Soc.}

\bibitem{B} J. Bourgain, {\it New classes of ${\mathscr
L}_p$-spaces,} Lecture Notes in Math. {\bf 889}, Springer, Berlin
1981.

\bibitem{Bo2} J. Bourgain, On the Dunford-Pettis property, {\it
Proc.\ Amer.\ Math.\ Soc.} {\bf 81} (1981), 265--272.

\bibitem{Bo3} J. Bourgain, Embedding $L^1$ in $L^1/H^1$, {\it
Trans.\ Amer.\ Math.\ Soc.} {\bf 278} (1983), 689--702.

\bibitem{BD} J. Bourgain and F. Delbaen, A class of special
${\mathscr L}_\infty$-spaces, {\it Acta Math.} {\bf 145} (1980),
155--176.

\bibitem{CG} J. M. F. Castillo and M. Gonz\'alez, On the
Dunford-Pettis property in Banach spaces, {\it Acta Univ.\
Carolin.--Math.\ Phys.} {\bf 35} (1994), 5--12.

\bibitem{DF} A. Defant and K. Floret, {\it Tensor Norms and
Operator Ideals,} Math.\ Studies {\bf 176}, North-Holland,
Amsterdam 1993.

\bibitem{DiDP} J. Diestel, A survey of results related to the
Dunford-Pettis property, in: W. H. Graves, {\it Proc.\ Conf.\ on
Integration, Topology and Geometry in Linear Spaces,} Chapel Hill
1979, Contemp.\ Math. {\bf 2}, 15--60, American Mathematical
Society, Providence RI 1980.

\bibitem{Di} J. Diestel, {\it Sequences and Series in Banach Spaces,}
Graduate Texts in Math.\ {\bf 92}, Springer, Berlin 1984.

\bibitem{DJT} J. Diestel, H. Jarchow and A. Tonge, {\it Absolutely
Summing Operators,} Cambridge Stud.\ Adv.\ Math. {\bf 43},
Cambridge University Press, Cambridge 1995.

\bibitem{DU} J. Diestel and J. J. Uhl, Jr., {\it Vector Measures,}
Math.\ Surveys Monographs {\bf 15}, American Mathematical Society,
Providence RI 1977.

\bibitem{Em1} G. Emmanuele, Remarks on weak compactness of
operators defined on certain injective tensor products, {\it
Proc.\ Amer.\ Math.\ Soc.} {\bf 116} (1992), 473--476.

\bibitem{Em2} G. Emmanuele, Some remarks on lifting of isomorphic
properties to injective and projective tensor products, {\it
Portugal.\ Math.} {\bf 53} (1996), 253--255.

\bibitem{Gr} A. Grothendieck, Sur les applications lin\'eaires
faiblement compactes d'espaces du type $C(K)$, {\it Canad.\ J.
Math.} {\bf 5} (1953), 129--173.

\bibitem{GrL1} A. Grothendieck, Une caract\'erisation
vectorielle-m\'etrique des espaces $L^1$, {\it Canad.\ J. Math.}
{\bf 7} (1955), 552--561.

\bibitem{GD} S. Guerre-Delabri\`ere, {\it Classical Sequences in
Banach Spaces,} Monographs and Textbooks in Pure and Appl.\ Math.
{\bf 166}, Dekker, New York 1992.

\bibitem{JR} W. B. Johnson and H. P. Rosenthal, On
$\omega^*$-basic sequences and their applications to the study of
Banach spaces, {\it Studia Math.} {\bf 43} (1972), 77--92.

\bibitem{K} N. J. Kalton, Spaces of compact operators, {\it Math.\
Ann.} {\bf 208} (1974), 267--278.

\bibitem{LR} J. Lindenstrauss and H. P. Rosenthal, The ${\mathscr
L}_p$-spaces, {\it Israel J. Math.} {\bf 7} (1969), 325--349.

\bibitem{LP} F. Lust, Produits tensoriels injectifs d'espaces de
Sidon, {\it Colloq.\ Math.} {\bf 32} (1975), 285--289.

\bibitem{OP} R. I. Ovsepian and A. Pe\l czy\'nski, On the
existence of a fundamental total and bounded biorthogonal sequence
in every separable Banach space, and related constructions of
uniformly bounded orthonormal systems in $L^2$, {\it Studia Math.}
{\bf 54} (1975), 149--159.

\bibitem{Pi} G. Pisier, {\it Factorization of Linear Operators and
Geometry of Banach Spaces,} Regional Conf.\ Ser.\ in Math. {\bf
60}, American Mathematical Society, Providence RI 1986.

\bibitem{Ra} N. Randrianantoanina, Some remarks on the
Dunford-Pettis property, {\it Rocky Mountain J. Math.} {\bf 27}
(1997), 1199--1213.

\bibitem{Ry} R. A. Ryan, The Dunford-Pettis property and
projective tensor products, {\it Bull.\ Polish Acad.\ Sci.\ Math.}
{\bf 35} (1987), 785--792.

\bibitem{S} C. Stegall, Duals of certain spaces with the
Dunford-Pettis property, {\it Notices Amer.\ Math.\ Soc.} {\bf 19}
(1972), A-799.

\bibitem{T} M. Talagrand, La propri\'et\'e de Dunford-Pettis dans
$C(K,E)$ et $L^1(E)$, {\it Israel J. Math.} {\bf 44} (1983),
317--321.

\end{thebibliography}
\end{document}